\documentclass{amsart}
\usepackage[latin1]{inputenc}
\usepackage[T1]{fontenc}
\usepackage[english]{babel}
\usepackage{amsmath,amssymb,epsf,epic,eepic,pifont}
\usepackage{graphicx}
\usepackage{amsfonts}
\usepackage[all]{xy}

\xyoption{arc}

\newtheorem{thm}{\bfseries Theorem}[section]
\newtheorem{lem}[thm]{\bfseries Lemma}
\newtheorem{prop}[thm]{\bfseries Proposition}

\newtheorem{defn}[thm]{\bfseries Definition}
\newtheorem{rem}[thm]{\bfseries Remark}

\renewcommand{\bf}{\bfseries}
\renewcommand{\it}{\itshape}

\numberwithin{equation}{section}

\newcommand{\D}{\mathcal{D}}

\newcommand{\s}{Cob}

\newcommand{\lc}{\Lambda_{\mathcal{C}}}

\newcommand{\oal}{\overline{\alpha}}

\newcommand{\p}{\Pi_1}
\newcommand{\z}{\mathcal{Z}}

\newcommand{\Z}{\mathbb{Z}}

\newcommand{\N}{\mathbb{N}}
\newcommand\acl{\left\{ \begin{array}{ccc} }
\newcommand\acr{\end{array}\right\} }
\newcommand{\m}{\mathcal}
\newcommand{\cbt}{Cob_{tri}}
\newcommand{\cb}{Cob}

\title{The invariant of Turaev-Viro from Group category}
\author{Jérôme Petit}
\address{Institut de Mathématiques et de Modélisation de Montpellier-U.M.R CNRS 5149
\newline Département des Sciences Mathématiques, Université de Montpellier~II,\newline Case
Courrier 051, Place Eugène Bataillon, 34095 Montpellier Cedex~5,
France.} \email{petit@math.univ-montp2.fr}
\email[url]{www.math.univ-montp2.fr/$\sim$petit/} 

\begin{document}

%%abstract---------------------
\begin{abstract}
A Group category is a spherical category whose simple objects are
invertible. The invariant of Turaev-Viro with this particular
category is in fact the invariant of Dijkgraaf-Witten whose the
group and the 3-cocycle is given by the simple objects and the
associativity constraint of the category.
\end{abstract}

\maketitle \tableofcontents
\section*{Introduction}
In 1992 M. Wakui \cite{Wakui} reformulated the invariant of
Dijkgraaf-Witten \cite{DW} and he proved the topological invariance
in a rigorous way. The invariance is based upon the triangulation
and the Pachner moves. Once given a finite group and a 3-cocycle the
Dijkgraaf-Witten invariant is defined combinatorially. Moreover in
this paper he built a topological quantum field theory (TQFT) from
this invariant. The same year V. Turaev and O. Viro \cite{TV} built
an invariant of 3 manifold thanks to 6-j symbol  to prove the
topological invariance they showed a relative version of a theorem
of Alexander \cite{Alex}  on equivalence of triangulation. This
invariant was reformulated in a categorical languages \cite{Tu} and
the TQFT was built. In the same spirit of \cite{Tu} J.W. Barret and
B.W. Westburry \cite{BW} have built a 3-manifold invariant using
spherical categories. In this construction the topological
invariance puts back down the triangulation and the Pachner moves.
Independently I. Gelfand and D. Kazhdan \cite{GK} have built a
3-manifold using spherical categories and in 1993 D.N. Yetter has
studied an untwisted version of the invariant of Dijkgraaf-Witten in
\cite{Yetter} the Turaev-Viro style in . In fact these constructions
are reformulations of the Turaev-Viro invariant. In the rest of the
paper we will call such kind of invariant the invariant of
Turaev-Viro and it will be denoted : $TV_{\m{C}}$ where $\m{C}$ is
the category used to build
the invariant.\\
The main goal of this paper is to give a relation between this two
approaches based on triangulation. That's why we utilize a "special"
spherical category. Roughly speaking, it is a spherical category
such that every simple object is invertible and has a dimension
equal to one. The dimension is given by the spherical structure. In
\cite{Quinn}, F. Quinn called this category : "Group category". In
\cite{FRS} invertible objects are called simples currents and the
tensor category whose every simple objects are invertible is denoted
Pointed category. The authors have denoted Picard category of
$\m{C}$ the full tensor category of $\m{C}$ whose objects are direct
sum of invertible objects of $\m{C}$. Thus if there is a finite
number of simple object and if every object is finite direct sum of
simple object then a pointed category is equal to its Picard
category. In this paper we will use the terminology of F. Quinn
\cite{Quinn}. L. Crane and D.N. Yetter have studied group cocycle to
describe monoidal category with duals in \cite{cryecat}.\\
 Here is an outline of the paper. In Section
\ref{rappelDW} we recall the definition of the Dijkgraaf-Invariant
\cite{Wakui}. In Section \ref{Group category} we give the definition
and we recall some facts on the Group category. In Section \ref{ITV}
we give the definition of the Turaev-Viro invariant of 3-manifold.
In Section \ref{egalite}, we compute the Turaev-Viro invariant in
the case of Group-category with other conditions and we show the
main theorem (\ref{eq})
. In section \ref{topcol} we give a topological interpretation of the admissible colorings. In Section \ref{TQFTTV} we give the construction the TQFT which arises from this invariant. We end the paper by discussing a few examples.\\
\section{The invariant of Dijkgraaf-Witten}\label{rappelDW}
Throughout this paper $k$ will be a commutative field such that
$car(k)=0$ and
$\overline{k}=k$.\\
We use the description of \cite{Wakui}. Let $G$ be a finite group,
this group will be always a multiplicative group. Moreover $k$ is a
representation of $G$ with the trivial action. Then we can define
$Z^3(G,k^{\star})$ the set of 3-cocycle of $G$ with coefficients in
$k^{\star}$ and we fix $\alpha \in Z^3(G,k^{\star})$. Let $T$ be a
$n$-simplex with $n\geq 1$, a color of $T$ is the following data :
\begin{equation}\label{color}
\gamma \mbox{ : }\{ \mbox{oriented edges of T } \} \rightarrow G,
\end{equation}
which satisfies the conditions :
\begin{itemize}
\item[(i)]for any oriented edge $e$ : $\gamma(\overline{e})=\gamma(e)^{-1}$, where
$\overline{e}$ is the oriented
edge with the opposite orientation. \\
\item[(ii)]For any oriented 2-simplex (012) of $T$ we have :
$$\gamma(01)\gamma(12)\gamma(20)=1.$$
\end{itemize}
We denote $Col(T)$ the set of all colors of $T$, if $T$ is a
triangulation of a n-manifold  $M$, with $n \geq 2$, we denote
$Col(M,T)$ the set of all colors of $M$ given by $T$. When there is
no ambiguity on the choice of a triangulation, we denote $Col(M)$
the set of colors of M. If $M$ is a manifold with boundary :
$\partial M$, then $\partial M$ is endowed with a triangulation
which comes from the triangulation of $M$. If $\tau$ is a color of
$\partial M$ then the set of all colors of $M$ which extend $\tau$,
is denoted $Col(M,\tau)$. We give an order to the set of vertices of
a triangulation of $M$, then each 3-simplex has an orientation given
by the ascending order. Then for $\gamma\in Col(M)$ and for the
3-simplex $(0123)$ we put :
\[\alpha(\Delta,\gamma)=\alpha(\gamma(01),\gamma(12),\gamma(20)),\]\label{alpha} with $\alpha \in Z^3(G,k^{\star})$.
\begin{thm}[Wakui (92)\cite{Wakui}]
Let $G$ be a finite group, we fix a 3-cocycle $\alpha \in
Z^3(G,k^{\star})$. Let $M$ be a compact oriented triangulated
3-manifold, $T$ is a triangulation of $M$. We denote the number of
vertices of $T$ by $n_0$ and $T^3$ the set of 3-simplex in $T$. All
the 3-simplex are oriented by a numbering of the vertices. Given
$\tau \in Col(\partial M)$, we define the Dijkgraaf-Witten invariant
by :
$$Z_M(\tau)=\mid G \mid ^{-n_0}\sum_{\gamma \in Col(M,\tau)}\prod_{\Delta \in
T^3}\alpha(\Delta, \gamma)^{\epsilon_{\Delta}},$$ where
$$
\epsilon_{\Delta}=\left\{
\begin{array}{cc}
1 & if \Delta \mbox{ and M have the same orientation}, \\
-1 & otherwise.
\end{array}
\right.
$$
Then $Z_M(\tau)$ does not depend on the choice of triangulation of
$M$ and the choice of order of vertices in $M$ whenever we fix a
triangulation of $\partial M$ and $\tau$.
\end{thm}
Thanks to the independence of the choice of numbering, we can
consider a numbering of the triangulation such that the 3-simplex
have the same orientation of $M$. Then the invariant is :
\[Z_M(\tau)=\mid G \mid ^{-n_0} \sum_{\gamma \in
Col(M,\tau)}\prod_{\Delta \in T^3}\alpha(\Delta, \gamma),\] where
all 3-simplex $\Delta$ have the same orientation of $M$.
\begin{rem}
If we consider $M$ without boundary, then $Z_M(\emptyset)$ is a
3-manifold invariant and we denote it : $Z_M$.
\end{rem}
\section{Group category}\label{Group category}
In this section, we review some basics facts on Group category.
\subsection{Definition}
Let $\m{C}$ be a monoidal category, by a scalar object
\cite{vir_kirby} of $\m{C}$ we shall mean an object of $\m{C}$ such
that : $End(X)=k$. If $\m{C}$ is abelian and $k$ is algebraically
closed then an object is scalar iff it is an simple object. We
denote the set of isomorphism classes of scalar objects of $\m{C}$
by $\lc$.

\begin{defn}
A finitely semisimple monoidal category is a monoidal category
$(\m{C},\otimes ,I,a,l,r)$ such that :
\begin{itemize}
\item[(a)] $\m{C}$ is an abelian $k$-category and $\otimes $ is a bifunctor $k$-linear,\\
\item[(b)] every object of $\m{C}$ is a finite direct sum of scalar objects of $\m{C}$,\\
\item[(c)] $ \sharp \lc <\infty$ and $I$ is a scalar object, \\
\item[(d)] $\m{C}$ is sovereign\footnote{$\m{C}$ admits a right and a left duality which
are isomorphic as monoidal functor.}
\end{itemize}
\end{defn}
If $\m{C}$ is a finitely semisimple monoidal category, then every
object $X$ of $\m{C}$ admits a right duality :
$(X,X^{\vee},e_X,h_X)$ and a left duality :
$(X^{\vee},X,\epsilon_X,\eta_X)$, we can take the same object
because $\m{C}$ is sovereign. By definition of duality :
\begin{align*}
e_X \mbox{ : }X\otimes X^{\vee} & \rightarrow I \\
\epsilon_X \mbox{ : }X^{\vee} \otimes X & \rightarrow I \\
\eta_X \mbox{ : } I & \rightarrow X\otimes X^{\vee} \\
h_X \mbox{ : }I & \rightarrow X^{\vee} \otimes X.
\end{align*}
and we have the following equalities :
\begin{align*}
(e_X \otimes id_X)(id_X \otimes h_X) & = id_X \\
(id_X \otimes \epsilon_X)(\eta_X \otimes id_X) & =id_X\\
(id_{X^{\vee}} \otimes e_X)(h_X \otimes id_{X^{\vee}})&=id_{X^{\vee}} \\
(\epsilon_X \otimes id_{X^{\vee}})(id_{X^{\vee}} \otimes \eta_X)&=
id_{X^{\vee}}.
\end{align*}
The left quantum trace of an endomorphism $f\in End_{\m{C}}(X)$ is
defined by :
$$tr_l(f)=e_X(f\otimes id_{X^{\vee}})\eta_X,$$
the right quantum trace of an endomorphism $f\in End_{\m{C}}(X)$ is
defined by :
$$tr_r(f)=\epsilon_X(id_{X^{\vee}}\otimes f)h_X.$$
for any endomorphisms $f,g$ in $\m{C}$ we have :
\begin{align*}
tr_l(f\otimes g)&=tr_l(f)tr_l(g), \\
tr_r(f\otimes g)&=tr_r(f)tr_r(g), \\
tr_r(f)&=tr_l(f^{\vee}),
\end{align*}
the multiplication is given by the multiplication of $k=End(I)$.
\begin{defn}\label{spherique}
A spherical category is a finitely semisimple monoidal category such
that, for all endomorphism $f$ in $\m{C}$ we have :
$tr_l(f)=tr_r(f)$.
\end{defn}
In a spherical category we denote the left trace by $tr$ and so we
have ${tr=tr_l=tr_r}$. The quantum dimension of an object $X$ in a
spherical category $\m{C}$ is defined by :
$$dim(X)=tr(id_X),$$
so we have $dim(X)=dim(X^{\vee})$.
\begin{defn}
\begin{itemize}
\item[]
\item[(i)]An object $X$ of a monoidal category $\m{C}$ is called invertible iff there exists
an object $Y$ such that $X\otimes Y \cong I$, where $I$ is the
tensor unit of $\m{C}$.
\item[(ii)]A monoidal category is called pointed iff every scalar object is
invertible.
\item[(iii)]The Group category {\it Pic($\m{C}$)} of the monoidal category $\m{C}$ is
the full monoidal subcategory of $\m{C}$ whose objects are direct
sums of invertible objects of $\m{C}$.
\item[(iv)]A Group category is a pointed spherical category.
\item[(v)] A $\theta$-category is a braided, pointed finitely semisimple monoidal category.
\end{itemize}
\end{defn}

\subsubsection{Example of Group category}

 $G$ is a finite group, we denote $k[G]$ the
category whose objects are $G$-graded finite dimensional $k$-vector
spaces\footnote{we can define a similar category, using $G$-graded
free $A$-modules, with $A$ a commutative ring} and whose morphisms
are $k$-linear morphism that preserves the grading. If $V$ and $W$
are objects of $k[G]$ the monoidal structure of $k[G]$ is given by :
 \[(V\otimes W)_g=\sum_{\begin{array}{c}
 h,k\\
 hk=g
\end{array}}V_h \otimes W_k.\] the associativity is
the identity and the isomorphism classes of scalar objects are in
bijection with $G$ : $g \leftrightarrow \delta_g$ where $\delta_g$
is defined in the following way :
$$
(\delta_g)_h= \left\{
\begin{array}{cc}
k & \mbox{if $g=h$}, \\
0 & otherwise.
\end{array}
\right.
$$
and every scalar object is invertible, thus $k[G]$ is a Group
category. $k[G]$ is a $\theta$-category iff $G$ is an abelian group.
\subsection{Some results on Group category}
Whenever $\m{C}$ is a Group category, it follows from the definition
of a Group category and the quantum dimension that for all $X \in
\lc$ : $dim(X)^2=1$. The Grothendieck ring of $\m{C}$ is isomorphic
to the group algebra of the finite group $\lc$, it is denoted
$\m{K}_0(\m{C})\cong \Z[\lc]$.
\begin{prop}
If $\m{C}$ is a Group category then :
\begin{itemize}
\item[(i)] all invertible objects are in $\lc$,\\
\item[(ii)] $(\lc,\otimes,I)$ is a finite group.
\end{itemize}
\end{prop}
{Proof (i) : }If $X$ is invertible in $\m{C}$ then there exists an
object $Y$ in $\m{C}$ such that : \newline${X\otimes Y \cong I}$,
thus we have :
$$\sum_{Z\in \lc}\mu_Z(X)Z\otimes Y=I,$$ where
$\mu_Z(X)=dim_k(Hom_{\m{C}}(X,Z))$ and so we have : \[\sum_{Z',Z\in
\lc}\mu_Z(X)\mu_{Z'}(Z\otimes Y)Z'=I,\] since $I$ is a scalar object
$$
\sum_{Z\in \lc}\mu_Z(X)\mu_{Z'}(Z\otimes Y)=\left\{
\begin{array}{cc}
1 & \mbox{if } Z'=I, \\
0 & otherwise.
\end{array}
\right.
$$
Moreover $Hom_{\m{C}}(X,Y)$ is finite for all objects in $\m{C}$ and
so $\mu_Z(X)\in \N \hookrightarrow k$, thus there is only one $Z_0
\in \lc$ such that $\mu_{Z_0}(X)\not=0$ and moreover
$\mu_{Z_0}(X)=1$ so $X=Z_0\in \lc$. We can notice that if $Y$ is the
inverse of $X$ then $X\otimes Y \cong I$ and so $Y\cong
X^{\vee}$.\qed

{Proof (ii) : } If $X$ is scalar then by definition of a Group
category $X$ is invertible and so there is $Y$ an object of $\m{C}$
such that $X\otimes Y \cong I$ we have seen that ${Y=X^{\vee} \in
\lc}$. In a finitely scalar monoidal category we have :
$X^{\vee}\otimes X=I \oplus Z$ where $Z$ is an object of $\m{C}$,
thus we have :
\begin{align*}
X^{\vee}& \cong X^{\vee} \otimes I \\
& \cong X^{\vee} \otimes (X \otimes X^{\vee}) \\
& \cong (X^{\vee} \otimes X) \otimes X^{\vee} \\
& \cong X^{\vee} \oplus Z\otimes X^{\vee}
\end{align*}
If $X$ is scalar then $X^{\vee}$ is scalar thus $Z\otimes
X^{\vee}=0$ and $X^{\vee}\not =0$. So it follows that $Z=0$ and
$X^{\vee}\otimes X \cong I$, then $X^{\vee}$ is the left and right
inverse of $X$. If $X$ and $Y$ are scalar objects then $X\otimes Y$
is an object of $\m{C}$ and :
$$End_{\m{C}}(X\otimes Y)\cong Hom_{\m{C}}(X,X\otimes Y\otimes Y^{\vee}) \cong End_{\m{C}}(X) \cong k,$$
then $X\otimes Y$ is a scalar object thus $(\lc,\otimes , I)$ is a
finite group.\qed

\begin{thm}[\cite{FK}, section 7.5]
Suppose $G$ is a finite group, then : Group categories with
underlying group $G$ correspond to $H^3(G,k^{\star})$.
\end{thm}
In fact $H^3(G,k^{\star})$ classifies all the associativity
constraint(up to monoidal equivalences). The group $G$ gives the set
of isomorphic classes of scalar objects and an element $\alpha \in
H^3(G,k^{\star})$ gives the associativity constraint of the Group
category. If we take $\alpha, \alpha' \in Z^3(G,k^{\star})$ such
that $[\alpha]=[\alpha'] \in H^3(G,k^{\star})$ then we obtain two
Group categories denoted by $\m{C}(G,\alpha)$ and $\m{C}(G,\alpha')$
such that : $\m{C}(G,\alpha)\cong^{\otimes}\m{C}(G,\alpha')$
(monoidal equivalence).
\subsection{6j-symbol}
We fix $\D$ a finitely monoidal category then for all object $X$ in
$\D$ we have : ${X=X_1\oplus ... \oplus X_n}$ with ${X_i \in
\Lambda_{\D}}$ then for all $1 \geq j \geq n$ there are morphisms
${i_j \in Hom_{\D}(X_j,X)}$ and ${p_j \in Hom_{\D}(X,X_j)}$ such
that ${p_ji_j=id_{X_j}}$ and ${\sum_ji_jp_j=id_X}$.
\begin{lem}\label{iso6j}
We fix $a,b,c,d,e,f \in \Lambda_{\D}$ then the following application
\begin{align*}
\Psi \mbox{ : } Hom(a,e\otimes d)\otimes_kHom(e,b\otimes c & \rightarrow Hom(a,(b\otimes c)\otimes d) \\
v\otimes w & \mapsto (w\otimes id_d)v
\end{align*}
induces an isomorphism between ${Hom_{\m{D}}(a,(b\otimes c)\otimes
d)}$ and
\newline ${\bigoplus_{e \in \Lambda} Hom_{\m{D}}(a,e\otimes d)\otimes_{k}
Hom_{\m{D}}(e,b\otimes c)}$. In the same vein we have :\newline
${Hom_{\m{D}}(a,b\otimes (c\otimes d)) \cong \bigoplus_{f \in
\Lambda} Hom_{\m{D}}(a,b\otimes f)\otimes_{k} Hom_{\m{D}}(f,c\otimes
d)}$
\end{lem}
{Proof :} By definition of $\D$ we have : $b\otimes c = \oplus_{e\in
\lambda_{\D}} \mu_e(b\otimes c)e$ with \\${\mu_e(b\otimes
c)=dim_k(Hom(e,b\otimes c))}$. Then for all $f\in
Hom_{\m{D}}(a,(b\otimes c)\otimes d)$ we have :
\begin{align*}
f &= id_{b\otimes c}\otimes id_d f \\
 &=\sum_{e \in \Lambda}(i_e p_e \otimes id_d)f \\
 &=\sum_{e \in \Lambda} (i_e\otimes id_d)(p_e\otimes id_d)f,
\end{align*}
and so $\Psi$ is surjective. Moreover the vector spaces are finite
dimensional and they have the same dimension thus we get the
isomorphism. The second isomorphism is obtained in the same way.\qed

$a$, the associativity constraint of $\D$, induces a natural
isomorphism : \newline ${(X\otimes Y)\otimes Z \cong X\otimes
(Y\otimes Z)}$, for all $X,Y,Z \in ob(\D)$. Then we have the
following commutative square :
$$
\xymatrix{\bigoplus_{e \in \Lambda} Hom_{\m{D}}(a,e\otimes d)\otimes_{k} Hom_{\m{D}}(e,b \otimes c)\ar[r] \ar[d]_{\cong}& \bigoplus_{f \in \Lambda} Hom_{\m{D}}(a,b\otimes f)\otimes_{k} Hom_{\m{D}}(f,c\otimes d)\ar[d]^{\cong}\\
Hom_{\m{D}}(a,(b\otimes c)\otimes d) \ar[r]_{\cong} &
Hom_{\m{D}}(a,b\otimes (c\otimes d))}
$$
the previous commutative square induces two linear applications :
\begin{align*}
&\acl a & b & c \\
     d & e & f \acr \mbox{ : }Hom_{\m{D}}(a,e\otimes d)\otimes_{k} Hom_{\m{D}}(e,b \otimes c) \rightarrow Hom_{\m{D}}(a,b\otimes f)\otimes_{k} Hom_{\m{D}}(f,c\otimes
     d) \\
&\acl a & b & c \\
     d & e & f \acr_{inv} \mbox{ : } Hom_{\m{D}}(a,b\otimes f)\otimes_{k} Hom_{\m{D}}(f,c\otimes
     d) \rightarrow Hom_{\m{D}}(a,e\otimes d)\otimes_{k} Hom_{\m{D}}(e,b \otimes c),
\end{align*}
these are the 6j-symbol of $\D$. \\
We define a bilinear form in the following way : for all objects
$X,Y$,
\begin{align*}
\omega_{X,Y} \mbox{ : }Hom_{\m{D}}(X,Y) \otimes Hom_{\m{D}}(Y,X) &\rightarrow k \\
f\otimes g & \mapsto tr_g(fg).
\end{align*}
By definition $\D$ doesn't admit
negligible morphism so $\omega_{\_,\_}$ is a non-degenerate bilinear form and it defines an adjoint of $\acl a & b & c \\
     d & e & f \acr$, for all $(a,b,c,d e,f)\in \Lambda_{\m{D}}$, this adjoint is denoted by : \[\lambda(a,b,c,d,e,f) \in (Hom_{\m{D}}(e\otimes d,a)\otimes
     Hom_{\m{D}}(b\otimes c,e)\otimes Hom_{\m{D}}(a,b\otimes f)\otimes Hom_{\m{D}}(f,c\otimes
     d))^{\star}.\]
\subsection{6j-symbol from Group category}
If $\m{C}$ is a Group category then for all $X, Y$ scalar objects
$X\otimes Y$ is a scalar object. Thus if $X,Y,Z$ are scalar objects
then :
\begin{equation}\label{evpicard}
Hom(Z,X\otimes Y)\cong \left\{
\begin{array}{cc}
k & \mbox{, if $X\otimes Y \cong Z$} \\
0 & ,otherwise
\end{array}
\right.
\end{equation}
In the case of Group category the isomorphisms (lemma \ref{iso6j})
become :
\begin{lem}\label{iso6jpic}

For all scalar objects $(a,b,c,d,e,f)$ we have :
\begin{align}
Hom_{\m{C}}(a,e\otimes d)\otimes_kHom_{\m{C}}(e,b\otimes c) & \cong Hom_{\m{C}}(a,(b\otimes c)\otimes d)\label{iso1} \\
Hom_{\m{C}}(a,b\otimes f)\otimes_kHom_{\m{C}}(f,c\otimes d )& \cong
Hom_{\m{C}}(a,b\otimes (c\otimes d))\label{iso2}
\end{align}
\begin{itemize}
\item[(i)]$\acl a & b & c \\d & e & f \acr \not= 0$ iff $e\cong b\otimes c$, $a \cong (b \otimes c)\otimes d$ and $f \cong c \otimes d$ \\
\item[(ii)]$\acl a & b & c \\d & e & f \acr_{inv} \not= 0$ iff $e\cong b\otimes c$, $a \cong b \otimes (c\otimes d)$ and $f \cong c \otimes d$
\end{itemize}
\end{lem}
{Proof : } The assertions (i), (ii) and the isomorphisms
(\ref{iso1}), (\ref{iso2}) come from (\ref{evpicard}). \qed
 Thus in the case
of the Group category the 6j-symbol $\acl a & b & c
\\d & e & f \acr$ only depends on $b,c,d$. For all scalar objects
$b,c,d$ we put ${\alpha(b,c,d)=\acl a & b & c \\d & e & f \acr}$. We
can define $\alpha $ in $\lc$. $\lc$ is a group and for $g \in \lc$
we denote $X_g$ a representation of this isomorphism class, and so
for all $g,h \in \lc$ we have : $X_g\otimes X_h \cong X_{gh}$. By
(\ref{evpicard}) we know that $Hom_{\m{C}}(X_{gh},X_g\otimes X_h)$
is a one dimensional vector space. For all $g,h \in \lc$ we put
$\phi(g,h)$ a basis of $Hom_{\m{C}}(X_{gh},X_g\otimes X_h)$.\\
We put $g,h,k\in \lc$ and we denote $X_g,X_h,X_k$ their
representations. By construction $\alpha(X_g,X_h,X_k)$ is an
isomorphism of one dimensional vector spaces thus in the basis
$\phi$ we have :
\begin{equation}
\alpha(X_g,X_h,X_k)(\phi(gh,k)\otimes
\phi(g,h))=\oal(g,h,k)(\phi(g,hk)\otimes \phi(h,k)),
\end{equation}
 with $\oal
\mbox{ : } \lc \times \lc \times \lc \rightarrow k^{\star}$. With
the same notations the commutative square which defines $\alpha$
induces the following equality :
\begin{equation*}
\oal(g,h,k)(id_{X_g}\otimes
\phi(h,k))\phi(g,hk)=a(X_g,X_h,X_h)(\phi(g,h)\otimes
id_{X_k})(\phi(gh,k)).
\end{equation*}
 Thus $\oal$ determine the following
isomorphism :
\begin{align*}
Hom_{\m{C}}(X_{ghk},(X_g \otimes X_h)\otimes X_k) & \cong Hom_{\m{C}}(X_{ghk},X_g \otimes (X_h\otimes X_k)) \\
v & \mapsto a(X_g,X_h,X_h)v,
\end{align*}
$a$ is the associativity constraint of $\m{C}$, thus $a$ satisfies
the Maclane's pentagon : with $X_g,X_h,X_k,X_l$ scalar objects
\[
\xy (0,0)*+{((X_g\otimes X_h)\otimes X_i)\otimes
X_j}="A";(25,20)*+{(X_g\otimes (X_h\otimes X_i))\otimes
X_j}="B";(50,00)*+{X_g\otimes(( X_h\otimes X_i)\otimes
X_j)}="C";(0,-20)*+{(X_g\otimes X_h)\otimes (X_i\otimes
X_j)}="D";(50,-20)*+{X_g\otimes (X_h\otimes (X_i\otimes X_j))}="E";
 {\ar^{a(g,h,i)\otimes id} "A";"B"};{\ar^{a(g,hi,j)} "B";"C"};{\ar^{id\otimes a(h,i,j)} "C";"E"};{\ar_{a(gh,i,j)}
"A";"D"};{\ar_{a(g,h,ij)} "D";"E"};

\endxy
\]

%$$a(X_g,X_h,X_k\otimes X_l)a(X_g \otimes X_h,X_k,X_l)=(id_{X_g}\otimes a(X_h,X_k,X_l))a(X_g,X_h\otimes
%X_k,X_l)(a(X_g,X_h,X_k)\otimes id_{X_l}).$$
If we apply the last equality in the basis $\phi$ we have :
$$\oal(g,h,kl)\oal(gh,k,l)=\oal(h,k,l)\oal(g,hk,l)\oal(g,h,k).$$

\begin{prop}
If $\m{C}$ is a Group category then the 6j-symbol is determined by a
3-cocycle on $Z^{3}(\lc,k^{\star})$ and a basis of
$Hom(X_{gh},X_g\otimes X_h)$.
\end{prop}
The relation of $a$ with the identity constraint $(r,l)$ induces
that : \newline ${l(h)\oal(g,1,h)=r(g)}$, for all $g,h \in \lc$ thus
$r(g)=\oal(g,1,1)$ and $\l(g)=\oal(1,1,h)^{-1}$. We can change
$\oal$ such that $\oal $ is normalized\footnote{r=l=1} without
changed the cohomologous class of $\oal$. In term of basis $\phi$,
it is a change of basis.
\section{The invariant of Turaev-Viro}\label{ITV}
We adopt the approach of \cite{GK} rather then \cite{TV}, but we use
a spherical category. Because if we consider a sovereign category
there is a problem in the construction. The problem occurs at the
level of independence of
the numbering of the 3-simplex.\\
Let $T$ be a $n$-simplex with $n\geq 1$, then a Turaev-Viro color of
$T$ (Turaev-Viro point of view) is the following data : $\gamma
\mbox{ : }\{ \mbox{oriented edges of T } \} \rightarrow \lc $ which
satisfies the conditions :
\begin{itemize}
\item[(i)]for any oriented edge $e$ : $\gamma(\overline{e})=\gamma(e)^{\vee}$, where
$\overline{e}$ is the oriented edge with the opposite orientation.
\end{itemize}
The set of all Turaev-Viro color of $T$ is denoted $Col_{TV}(T)$.
let $T$ be a $n$-simplex and we fix a numbering of the vertices of
$F$, every faces of $T$ has an orientation given by the ascending
order : $(012)$. For every faces $(012)$ we define the following
vector space :
$V_{\m{C}}((012),\gamma)=Hom_{\m{C}}(I,\gamma(01)\otimes \gamma(12)
\otimes \gamma(20))$.
\begin{lem}
\begin{align}
V_{\m{C}}((012),\gamma)&\cong V_{\m{C}}((201),\gamma) \cong V_{\m{C}}((120),\gamma) \label{pos}\\
V_{\m{C}}((012),\gamma)&\cong
V_{\m{C}}((021),\gamma)^{\star}\label{neg}
\end{align}
\end{lem}
{Proof (\ref{pos}) : } It comes from the sovereign structure of
$\m{C}$.\newline For all $X,Y,Z \in ob(\m{C})$ we have :
\begin{align*}
Hom_{\m{C}}(I,X\otimes Y \otimes Z) &\leftrightarrow Hom_{\m{C}}(I,Y\otimes Z \otimes X) \\
f & \mapsto (\epsilon_X \otimes id_{Y\otimes Z \otimes
X})(id_{X^{\vee}}\otimes f \otimes id_X)(h_x) \qed
\end{align*}

 {Proof (\ref{neg}) : }
It comes from the fact that the category $\m{C}$ doesn't admit
negligible morphism and so the following bilinear form is
non-degenerate :
\begin{align*}
contr \mbox{ : }V_{\m{C}}((012),\gamma) \otimes V_{\m{C}}((021),\gamma) & \rightarrow k \\
f \otimes g & \mapsto f^{\vee}g=tr(f^{\vee}g)=tr(g^{\vee}f)\qed
\end{align*}

Thus the vector space $V_{\m{C}}((012),\gamma)$ is independent of
the starting point and if we change the orientation of the 2-simplex
we obtain the dual vector space. Moreover this dual vector space can
be obtain by a change of color in fact :
$$V_{\m{C}}((012),\gamma) \cong V_{\m{C}}((021),\gamma)^{\star} \cong V_{\m{C}}((021),\gamma'),$$
with $\gamma'(02)=\gamma(01)$, $\gamma'(21)=\gamma(12)$,
$\gamma'(10)=\gamma(20)$. Let $T$ be the triangulation of a compact
oriented surface $\Sigma$ and $T^2$ the set of 2-simplex of $T$,
then we define \[ V_{\m{C}}(\Sigma,T)=\bigoplus_{\gamma \in
Col(T)}\bigotimes_{f \in T^2}V_{\m{C}}(f,\gamma),\]\label{Vsurf}
and this space is independent of the choice of a numbering of $T$. \\
Let $\Delta$ be a 3-simplex, a numbering of the vertices of $\Delta$
gives an orientation of $\Delta$, with this orientation $\Delta$ is
denoted $(0123)$. We take $\gamma \in Col(\Delta)$ and we put :
\begin{align}
& V_{\m{C}}((132),\gamma)\otimes V_{\m{C}}((023),\gamma)\otimes V_{\m{C}}((031),\gamma)\otimes V_{\m{C}}((012),\gamma)  \stackrel{L((0123),\gamma)}{\longrightarrow} k \label{L}\\
& v_0 \otimes v_1 \otimes v_2 \otimes v_3  \mapsto
dim(\gamma(13))^{-1}\lambda(\gamma(03),\gamma(01),\gamma(12),\gamma(23),\gamma(02),\gamma(13)).
\nonumber
\end{align}
This application defines, with duality given by $\omega$, an element
\newline  ${\tilde{L}((0123),\gamma)\in V_{\m{C}}((132),\gamma)^{\star}\otimes
V_{\m{C}}((023),\gamma)^{\star}\otimes
V_{\m{C}}((031),\gamma)^{\star}\otimes
V_{\m{C}}((012),\gamma)^{\star}}$. If $T$ is a triangulation of $M$,
which is an oriented and closed 3-manifold, then we denote $T^3$ the
set of oriented $3$-simplex of $T$ and we define the following
element :
\[\bigotimes_{\sigma \in T^3} \tilde{L}(\sigma, \gamma).\] But $M$ is a closed
3-manifold so every 2-simplex is a face of exactly two $3$-simplex
with opposite orientation. The 3-simplex are oriented such that
their orientations correspond to the orientation of $M$. We denote
$f$ the common face of $\sigma_1$ and $\sigma_2$, so the elements
can be written in the following way :  $\tilde{L}(\sigma_1,\gamma)
\in W \otimes V_{\m{C}}(f,\gamma)$, where $W$ is the tensor product
of the three other faces and $\tilde{L}(\sigma_2,\gamma) \in W'
\otimes V_{\m{C}}(\overline{f},\gamma)$ with $W'$ the tensor product
of the three other faces. The sovereign structure of $\m{C}$ defines
a bilinear non-degenerate form on this two vector spaces :
\begin{align*}
contr \mbox{ : } V_{\m{C}}((012),\gamma)\otimes V_{\m{C}}((021),\gamma) & \rightarrow k \\
f\otimes g & \mapsto tr(f^{\vee}g)=tr(g^{\vee}f),
\end{align*}
the  equality comes from the fact that $I=I^{\vee}$ and the
semi-simplicity of $\m{C}$ implies the non degeneracy of $contr$.
Since the 3-manifold is closed we can contract every 2-simplex and
then
\[Z(M,\gamma)=contr(\otimes_{\sigma \in T^3}
\tilde{L}(\sigma,\gamma)) \in k.\] We fix as in \ref{rappelDW},
$n_0$ the number of vertices of a given triangulation (we don't call
it $n_0(T)$ for two reasons, the first is historical \cite{TV},
\cite{Tu}, \cite{Wakui}, \cite{GK}, \cite{BW} and the second reason
comes from the fact that we use $n_0$ to describe an object which
doesn't depend on the triangulation). $T_o$ is the triangulation
with the orientation given by the numberings of the vertices such
that the orientaton is the orientation of the manifold $M$. The
invariant of Turaev-Viro is :
\begin{eqnarray}\label{invTV}
TV(M)=(\sum_{X\in \lc}dim(X)^2)^{-n_0}\sum_{\gamma \in
Col_{TV}(T)}\prod_{e\in T_o^1}dim(\gamma(e))Z(M,\gamma),
\end{eqnarray}
in the rest of the paper $\sum_{X\in \lc}dim(X)^2$ will be denoted
by $dim(\m{C})$.
\section{The equality}\label{egalite}
The invariant (\ref{invTV}), we make a sum over the Turaev-Viro
coloring of $T$ and we compute $L(\Delta,\gamma)$ for each
Turaev-Viro color $\gamma$ and each $3$-simplex $\Delta$. The linear
$L(\Delta,\gamma)$ is computed over vector spaces which are :
$V(f,\gamma)$ for each face $f$ of the 3-simplex $\Delta$. But in a
monoidal semisimple category we have the following result for all
scalar objects $a,b,c$ of $\m{C}$ :
$$
Hom_{\m{C}}(I,a\otimes b \otimes c)\left\{
\begin{array}{cc}
\ncong 0 &, a^{\vee}  \hookrightarrow b \otimes c \\
\cong 0 &, otherwise
\end{array} \right.
$$
If $\m{C}$ is Group category then we have the following relation :
$$
Hom_{\m{C}}(I,a\otimes b \otimes c)\left\{
\begin{array}{cc}
\cong k &, a^{\vee} \cong b \otimes c \\
\cong 0 &, otherwise
\end{array} \right.
$$
\begin{defn}\label{admi}
Let $\m{C}$ a Group category an $T$ an $n$-simplex with $n\geq1$. An
admissible colouring of $T$ is the set of application $\gamma$ from
oriented edges of $T$ to $\lc$ which satisfy :
\begin{itemize}
\item[(i)] $\gamma(\overline{e})=\gamma(e)^{\vee}$, where $\overline{e}$ is the edge $e$ with the opposite
orientation\\
\item[(ii)]for any oriented 2-simplex of $T$ we have :
$$\gamma(01)\otimes \gamma(12) \otimes \gamma(20) \cong I$$
\end{itemize}
\end{defn}
If $\m{C}$ is a Group category, then an admissible coloring is
nothing else than a color in a sense of Wakui (\ref{color}). That's
why we denote it $Col(T)$. If $\gamma \in Col_{TV}$ and $\gamma \not
\in Col(T)$, there is at least one oriented face $(012)$ in $T$ such
that : $\gamma(01)\gamma(12)\not= \gamma(02)$. It result that :
$V(012,\gamma)=0$. In a closed manifold every face is in the
boundary of exactly two 3-simplex with opposite orientation. the
value of $L(\_,\gamma)$ on these 3-simplex is equal to $0$. Thus in
the sum \ref{ITV} we have $Z(T,\gamma)=0$ for every $\gamma \not \in
Col(T)$. That's why can consider the sum \ref{ITV} only on the
admissible coloring (or coloring in Wakui sense).
\begin{thm}\label{eq}
Let $\m{C}$ a Group category such that for all $X \in \lc$ we have : $dim(X)=1$. \\
If $G$ is the underlying group of $\m{C}$ and if $\alpha \in
Z^3(G,k^{\star})$ is the associativity constraint of $\m{C}$ then
for all closed and oriented 3-manifold $M$ :
$$DW_{G,\alpha}(M)=TV_{\m{C}}(M).$$
\end{thm}
{Proof : } $\m{C}$ is Group category The condition $dim(X)=1$, for
all $X \in \lc$ implies : $dim(\m{C})=\sharp\lc$ and so :
$$TV_{\m{C}}(M)=(\sharp \lc)^{-n_0}\sum_{\gamma \in Col_{TV}(T)}Z(T,\gamma),$$
it remains to compute $Z(T,\gamma)$ for an admissible coloring.
\begin{lem}
If $\m{C}$ is a Group category such that for all $X \in \lc$
$dim(X)=1$ then : \[Z(T,\gamma)=\prod_{(0123) \in
T^3}\alpha(\gamma(01),\gamma(12),\gamma(23)),\] with $T$ a
triangulation of a closed and oriented 3-manifold and $\gamma$ an
admissible coloring of $T$.
\end{lem}
{Proof : } If $T$ is a triangulation of a closed and oriented
3-manifold and $\gamma$ is a coloring of $T$ then :
\[Z(T,\gamma)=contr(\otimes_{(0123)\in
T^3}\tilde{L}((0123),\gamma)\in k.\] Moreover if $dim(X)=1$ for all
$X \in \lc$ , by definition of $L$ (\ref{L}) we have
${L((0123),\gamma)=\lambda(\gamma(03),\gamma(01),\gamma(12),\gamma(23),\gamma(13),\gamma(02))}$.
$\m{C}$ is a Group category then $V((012),\gamma)\cong
Hom(\gamma(02),\gamma(01)\otimes \gamma(12))\cong k$. We fix $
\Phi(\gamma(01),\gamma(12))$ a basis of this vector space. If we
consider $(021)$, the same face with the opposite orientation, then
we have, thanks to the contraction : $V((012),\gamma)^{\star} \cong
V((021),\gamma)$ and so we can take the dual basis, it induces a
basis $\Phi'(\gamma(02),\gamma(21))$ of $Hom_{\m{C}}(\gamma(02)
\otimes \gamma(21), \gamma(01)) \cong V(021,\gamma)$
$$\acl a & b &c \\
d & e & f \acr(\Phi(b\otimes c,d)\otimes
\Phi(b,c))=\alpha(b,c,d)\Phi(b,c\otimes d)\otimes \Phi(c,d),$$ and
the duality given by the non-degenerate bilinear form $\omega$ gives
a basis of $Hom(e\otimes d, a)$ which is the dual of $Hom(a,e\otimes
d)$ :
\begin{align*}
\omega \mbox{ : } Hom(e \otimes d,a)\otimes Hom(a,e\otimes d) & \rightarrow k \\
f\otimes \Phi(e,d) & \mapsto tr(f\Phi(e,d))
\end{align*}
and so $\Phi'(e,d)=\Phi(e,d)^{-1}$. Thus we have :
\begin{align*}
&\tilde{L}((0123),\gamma) = \tilde{\lambda}(\gamma(03),\gamma(01),\gamma(12),\gamma(23),\gamma(13),\gamma(02)) \\
&=
\alpha(\gamma(01),\gamma(12),\gamma(23))\Phi(\gamma(01)\otimes\gamma(12),\gamma(23))^{-1}\otimes
\Phi(\gamma(01),\gamma(12))^{-1} \\
&\otimes \Phi(\gamma(01),\gamma(12)\otimes \gamma(23))\otimes
\Phi(\gamma(12),\gamma(23))
\end{align*}
\cite{GK} asserts that $L(0123,\gamma)$ doesn't depend on the choice
of the numbering which preserve the orientation of $(0123)$
\footnote{In \cite{GK} the authors asserts this result for a
sovereign category, but there is a problem with this condition. If
we use spherical category the result is true.} so we can choose any
numbering (0123) which preserve the orientation of the 3-simplex.
The contraction on $V((012),\gamma)$ induces a contraction on
$Hom(\gamma(01)\otimes \gamma(12),\gamma(02))$ given by the
isomorphism : $V((012),\gamma)\cong Hom(\gamma(01)\otimes
\gamma(12),\gamma(02))$ :
\begin{align*}
contr \mbox{ : }Hom(\gamma(01)\otimes \gamma(12),\gamma(02))\otimes
Hom(\gamma(02),\gamma(01)\otimes \gamma(12)) &
\rightarrow k \\
f \otimes g &\mapsto tr (fg)
\end{align*}
If we consider two 3-simplexes $(0123)$ and $(0214)$ with the same
orientation, they have a common face. We change the numbering of
$(0214)$ without changing the orientation : $(4012)$ (the cycle
$(042)$ has an even signature).
\begin{align*}
&\tilde{L}((0123),\gamma) = \tilde{\lambda}(\gamma(03),\gamma(01),\gamma(12),\gamma(23),\gamma(13),\gamma(02)) \\
&=
\alpha(\gamma(01),\gamma(12),\gamma(23))\Phi(\gamma(01)\otimes\gamma(12),\gamma(23))^{-1}\otimes
\Phi(\gamma(01),\gamma(12))^{-1} \\
&\otimes \Phi(\gamma(01),\gamma(12)\otimes \gamma(23))\otimes
\Phi(\gamma(12),\gamma(23))\\
&\tilde{L}((4012),\gamma) = \tilde{\lambda}(\gamma(42),\gamma(40),\gamma(01),\gamma(12),\gamma(02),\gamma(41)) \\
&=\alpha(\gamma(40),\gamma(01),\gamma(12))\Phi(\gamma(40)\otimes\gamma(01),\gamma(12))^{-1}\otimes
\Phi(\gamma(40),\gamma(01))^{-1} \\
&\otimes \Phi(\gamma(40),\gamma(01)\otimes \gamma(12))\otimes
\Phi(\gamma(01),\gamma(12))
\end{align*}
\begin{align*}
contr(\Phi(\gamma(01),\gamma(12))^{-1}\otimes
\Phi(\gamma(01),\gamma(12))=dim(\gamma(01)\otimes \gamma(12))=1
\end{align*}
$Z(T,\gamma)=\prod_{(0123) \in
T^3}\alpha(\gamma(01),\gamma(12),\gamma(23))$.\qed

\section{Topological interpretation of admissible
coloring}\label{topcol} In this section, we will give a topological
interpretation of the admissible coloring of a triangulation.
\subsection{The fundamental groupoïd of $T$}
Let $T$ be a n-simplex, we denote $\p(T)$ the following category :
\begin{itemize}
\item[]$Ob(\p(T))=T^0$
\item[]Arrows of $\p(T)$ are the oriented edges of $T$ and the $0$-simplex modulo the relation of
2-simplex, that is if $(012)$ is an oriented 2-simplex then
$(01).(12)=(02)$.
\end{itemize}
The composition is given by the concatenation of the edge, and the
inverse of an edge is the same edge with the opposite orientation.
The identity is given by the $0$-simplex himself. We can define the
pointed fundamental groupoïd $\p(T,x)$ in the same way. There is
only one object which is $x$ and the set of arrows are loops in x
and a loop is a concatenation of edges.
\begin{rem}
If $T$ is the triangulation of a connected manifold $M$, then there
is an equivalence of category between $\p(T)$ and the pointed
category $\p(T,x)$ where $x$ is $0$-simplex of $T$. Moreover the set
of arrows of $\p(T,x)$ is $\p(M,x)$.
\end{rem}
If $G$ is a group then the groupoïd obtained thanks to $G$ will be denoted $\m{G}$.\\
If $\m{C}$ is a Group category we can define the following
application :
\begin{align*}
\Psi \mbox{ : } Col(T) & \rightarrow Fun(\p(T),\lc) \\
\gamma & \mapsto F_{\gamma},
\end{align*}
the functor $F_{\gamma}$ is defined in the following way : for all
$x \in T^0$ we have ${F_{\gamma}(x)=\star}$, which is the object of
$\m{G}$. And $F_{\gamma}(01)=\gamma(01)$. $\gamma$ respects the
2-simplex condition, thus $F_{\gamma}$ is well defined.
\begin{lem}
$\Psi$ is bijective.
\end{lem}
{Proof : } If $F_{\gamma}=F_{\theta}$ then for all oriented edges
$e$ we have :
$\gamma(e)=\theta(e)$ and so $\gamma=\theta$. \\
If $F \in Fun(\p(T),\lc)$ then for all oriented edge $e$ we have
$F(e)\in \lc$ and $F(\overline{e})=F(e)^{\vee}$. Thus we can define
a coloring of $T$ : $\gamma(e)=F(e)$. We have to check the 2-simplex
condition. If $(012)$ is a 2-simplex then :
${F((02))=F((01).(12))=F((01))\otimes F((12))}$, thus $\gamma \in
Col(T)$ and $\Psi(\gamma)=F_{\gamma}$ and for object $F_{\gamma}=F$
and for all arrow $e$ : $F_{\gamma}(e)=\gamma(e)=F(e)$ thus $\Psi$
is bijective. \qed

\subsection{The gauge action}\label{gauge}
We note $\lc^{T^0}$ the set of application from the 0-simplex of $T$
to $\lc$. We can define an action of $\lc^{T^0}$ on $Col(T)$ in the
following way :
\begin{align*}
\lc^{T^0} \times Col(T) & \rightarrow Col(T) \\
(\delta,\gamma) & \mapsto \gamma^{\delta},
\end{align*}
such that for all oriented edge $(01)$ :
${\gamma^{\delta}(01)=\delta(0)\otimes \gamma(01)\otimes
\delta(1)^{\vee}}$. By a straightforward computation we
show that it is an action. \\
We define the following equivalence relation on $Col(T)$ :
$$(\gamma \sim \gamma') \Leftrightarrow (\exists \delta \in \lc^{T^0} \mbox{ such that }
\gamma^{\delta}=\gamma'),$$ we can construct the following
application :
\begin{align}\label{theta}
\Theta\mbox{ : } \frac{Col(T)}{\sim} &\rightarrow \frac{Fun(\p(T),\lc)}{iso} \\
[\gamma] & \mapsto [\psi(\gamma)]
\end{align}
This application is well defined because if
$\gamma'=\gamma^{\delta}$ then : ${\beta(x)=\delta(x)^{\vee}\mbox{ :
} F_{\gamma}(x) \rightarrow F_{\gamma'}(x)}$ is an in isomorphism in
the groupoïd $\lc$. For all oriented edge $(01)$ we have :
\begin{align*}
\beta(1)F_{\gamma}(01) &=  \gamma(01)\otimes \delta(1)^{\vee} \\
&= \delta(0)^{\vee}\otimes \delta(0) \otimes \gamma(01) \otimes \delta(1)^{\vee} \\
&=\delta(0)^{\vee}\gamma^{\delta}(01) \\
&=F_{\gamma^{\delta}}(01)\beta(0).
\end{align*}
Thus $\beta$ is a natural isomorphism between $F_{\gamma}$ and
$F_{\gamma^{\delta}}$.
\begin{prop}
$\Theta$ is a bijection.
\end{prop}
{Proof : }
$\Psi$ is surjective thus it follows that $\Theta$ is surjective. \\
Let $\gamma$ and $\gamma'$ two admissible colorings of $T$, if
$\Theta(\gamma)=\Theta(\gamma')$ then there is a natural isomorphism
between $F_{\gamma}$ and $F_{\gamma'}$. We note $\beta$ this
isomorphism, for all $0$-simplex $x$ we have : $\beta(x) \in \lc$
and $\beta(x) \mbox{ : } F_{\gamma}(x)\cong F_{\gamma'}(x)$. So
$\beta \in \lc^{T^0}$ and for all oriented edges $(xy)$ we have the
following commutative square :
$$
\xymatrix{ F_{\gamma}(x) \ar[r]^{\gamma(xy)}\ar[d]^{\beta(x)} & F_{\gamma}(y) \ar[d]^{\beta(y)} \\
F_{\gamma'(x)} \ar[r]_{\gamma'(xy)} & F_{\gamma'(y)} }
$$
thus $\gamma'^{\beta}=\gamma$. \qed

\section{Construction of TQFT}\label{TQFTTV}
\subsection{Triangulate TQFT}
We recall some definitions on TQFT and some results of \cite{Petit}.
We denote $\cb$ the category of 1+2 cobordism : $Ob(\s)$ is the set
of oriented and closed surface, the morphism of $\cb$ are the class
of oriented compact 3-manifold, i.e : the classes of diffeomorphisms
preserving the boundary. The disjoint union and $\emptyset $ give a
strict monoidal structure to $\cb$.
\begin{defn}\label{TQFT}
 $A$ is a commutative ring with unit, a TQFT is a monoidal and $A$-linear functor from $\cb$ to $A$-mod.
\end{defn}
In \cite{Petit}, there is a way of obtaining a TQFT from a functor
which is not monoidal.
\begin{defn}\label{unitaire}
Let $\m{C}$ a monoidal category, a monoidal, non unitary functor $F$
is the data\newline $(F,\Phi_2,\Phi_0) \mbox{ : } \mathcal{C}
\rightarrow \mathcal{D}$. $(F,\Phi_2,\Phi_0)$ verify all the axioms
of a monoidal functor expected the following : there is at least one
$X\in Ob(\mathcal{C})$ such that : $F(id_X)\not= id_{F(X)}$
\end{defn}

\begin{prop}[\cite{Petit}]\label{unitaire}
Let $(F,\Phi_2,\Phi_0)\mbox{ : }\m{C}\rightarrow \D$ a monoidal non
unitary functor, with $\mathcal{D}$ a monoidal, $A$-linear and
abelian category, then there is a monoidal functor $\tilde{F}\mbox{
: }\m{C}
 \rightarrow \D$ such as : for all $X\in Ob(\m{C})$ $\tilde{F}(X)$
is an sub object of $F(X)$.
\end{prop}
\subsection{Construction of Turaev-Viro}
We recall the construction of Turaev-Viro \cite{TV} and we apply the result of \cite{Petit}.\\
\underline{First step} \\
For every object $(\Sigma,T)$ we assign finite vector space :
$$V(\Sigma,T)=\bigoplus_{c \in Col(T)}\bigotimes_{f \in T^2}V(f,\gamma)=k[Col(T)],$$
it is the vector space spanned by the admissible coloration of $T$.\\
For every 3-manifold $M$  whose boundary is
$(-\Sigma,T)\coprod(\Sigma',T')$ and for every admissible coloring
c,c' of $T$ and $T'$ : $TV_M(c,c')\in k$. Thus we can define :
\begin{align}
V(M) \mbox{ : } V(\Sigma,T) & \rightarrow V(\Sigma',T')\label{V} \\
c & \mapsto \sum_{c' \in Col(T')}TV_M(c,c')c'\nonumber
\end{align}
By construction :
\begin{align}
V(M)V(N)(c)&=\sum_{c_1,c_2}TV_N(c,c_1)TV_M(c_1,c_2)c_2 \nonumber\\
&= (\sharp\lc)^{n_0(\partial N_+)}\sum_{c_2}TV_{M\circ N}(c,c_2)c_2\nonumber \\
&=(\sharp\lc)^{n_0(\partial N_+)}V(M\circ N)(c).\label{anomaly}
\end{align}
There are at least three ways of erasing the anomaly. Here are the
normalization, with $M$ a 3-manifold whose the boundary is $\partial
M=\overline{M}_-\coprod M_+$ and $c\in Col(M_-),\mbox{ }c'\in
Col(M_+)$ .
\begin{align}
TV_i(M)(c,c')&=\lc^{-n_0(M_-)}TV(M)(c,c') \\
TV_o(M)(c,c')&=\lc^{-n_0(M_+)}TV(M)(c,c') \\
TV_m(M)(c,c')&=\lc^{\frac{-n_0(M_-)-n_0(M_+)}{2}}TV(M)(c,c')
\end{align}
\begin{lem}\label{normalization}
\begin{itemize}
\item[]
\item[(i)] $TV_i$, $TV_o$ and $TV_m$ are invariants of 3-manifold with
boundary.
\item[(ii)] $TV_i$,$TV_o$ and $TV_m$ define the same monoidal non unitary functor (up to monoidal equivalence).
\end{itemize}
\end{lem}
{Proof (i) : }The theorem of Pachner define invariant of 3-manifold
whose the boundary is fix, thus the triangulation of the boundary
remains unchanged. That's why we obtain an invariant of 3 -manifold
with boundary.

{Proof (ii) : } Let $V_i$ (resp. $V_o$, $V_m$) the non unitary
functor defines from the invariant $TV_i$ (resp. $TV_o$, $TV_m$).
The natural transformation :
\begin{align*}
\beta \mbox{ : }V_o(\Sigma,T)=k[Col(T)]&\rightarrow V_i(\Sigma,T)=k[Col(T)] \\
c & \mapsto (\sharp\lc)^{n_0((\Sigma,T))}c
\end{align*}
is an isomorphism. It remains to show that $\beta$ is monoidal.
Since $V_o$ and $V_i$ are strict we the following square :
$$
\xymatrix{V_o((\Sigma,T)\otimes(\Sigma',T'))\ar[d]_{=}\ar[rrr]^{\beta((\Sigma,T)\otimes(\Sigma',T'))}&
& & V_i((\Sigma,T)\otimes(\Sigma',T'))\ar[d]^{=}
\\
V_o(\Sigma_1,T_1)\otimes V_o(\Sigma_1',T_1')
\ar[rrr]^{\beta(\Sigma,T)\otimes \beta(\Sigma',T')}& & &
V_i(\Sigma,T)\otimes V_i(\Sigma',T') }
$$
and for all $c\in Col(T\coprod T')=Col(T)\otimes Col(T')$,
\begin{align*}
\beta((\Sigma,T)\otimes(\Sigma',T'))(c)&=(\sharp
\lc)^{n_0((\Sigma,T)\coprod (\Sigma',T'))}c
\\
&=(\sharp \lc)^{n_0((\Sigma,T))+n_0((\Sigma',T'))}c \\
&=(\sharp \lc)^{n_0((\Sigma,T))+n_0((\Sigma',T'))}c_1\otimes c_2,
\end{align*}
with $c=c_1\otimes c_2$. And for $c\in Col(T), c'\in Col(T')$, we
have :
$$
\beta((\Sigma,T))\otimes\beta((\Sigma',T'))(c\otimes c')=(\sharp
\lc)^{n_0((\Sigma,T))}(\sharp \lc)^{n_0((\Sigma',T'))}c\otimes c'.
$$
Thus the square
commutes.\\
We prove in the way the monoidal isomorphism between $V_o$ and
$V_m$, the isomorphism is given by :
\begin{align*}
\kappa((\Sigma,T)) \mbox{ : } V(\Sigma,T) & \rightarrow V(\Sigma,T)
\\
c &\mapsto (\sharp \lc)^{-\frac{n_0(\Sigma,T)}{2}}c
\end{align*}
\qed  In \cite{Wakui} and \cite{Tu}, the authors used $TV_m$, here
we will $TV_i$ because in this case we don't have to compute the
square of $Dim(\m{C})$. We replace $TV$ by $TV_i$ in the definition
of (\ref{V}), and we still denote it by $V$.
\begin{prop}[\cite{Tu}, \cite{Wakui}]
$V$ is a monoidal non unitary functor from $\cbt$ to the category of
finite dimensional vector spaces.
\end{prop}

We fix $V(\Sigma\times I,T)$ the linear application given by the
identity on $(\Sigma,T)$ in $\cbt$.  \cite{Wakui} and \cite{Tu}, or
in a general framework \cite{Petit}, the TQFT is the following
functor :
\begin{align*}
\m{V} \s & \rightarrow k-vect \\
\Sigma & \mapsto \m{V}(\Sigma)=im(V(\Sigma\times I,T)) \\
M\in Hom(\Sigma,\Sigma') & \mapsto \m{V}(M)=V(M)_{\mid
im(V(\Sigma\times I,T))}
\end{align*}
The functor is well defined because for all triangulations $T$ and
$T'$ of $\Sigma$ we have : $im(V(\Sigma\times I,T))\cong
im(V(\Sigma\times I,T'))$. The isomorphism is given by $\Sigma\times
I$ with $T$ the triangulation of $\Sigma\times \{ 0 \}$ and $T'$ the
triangulation $\Sigma \times \{ 1\}$.
\section{Examples}\label{calcul}
\subsection{$\alpha=1$}
If $\alpha=1$ then $\m{C}$ is a strict Group category and $\m{C}$ is
equivalent to $k[G]$. Thanks to the gauge action (\ref{gauge}) and
the isomorphism (\ref{theta}), we have for all closed, oriented and
connected 3-manifold $M$ :
$$TV_{k[G]}(M)=\#G^{-n_0}\#Col(T)=\#\frac{Col(T)}{\sim}=\#\frac{Fun(\p(T,x),\lc)}{iso},$$
the groupoïd has only one object and so a functor is only defined by
the applications from the $\p(M,x)$ to $\lc$ and an isomorphism
between two functors implies that the applications are conjugate. So
we have :
$$TV_{k[G]}(M)=\#\frac{Hom(\p(M),\lc)}{conj}.$$
\subsection{$G=\Z_n$}
If $G$ is the cyclic group of order $n$ then we have :
$H^3(G,k^{\star})=\Z_n$, and $\alpha$ given by (\ref{zn})
is a 3-cocycle. \\
A triangulation of $S^1\times S^1\times S^1$ is given in \cite{DW} :

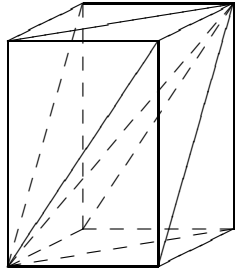
\[
\xy (0,0)*{}="A"; (20,0)*{}="B";
(0,30)*{}="C";(20,30)*{}="D";(10,35)*{}="E";(30,35)*{}="F";(30,05)*{}="G";(10,05)*{}="H";
 "A"; "B" **\dir{-}; "A";
"C" **\dir{-}; "C"; "D" **\dir{-}; "B";"D" **\dir{-}; "C";"E"
**\dir{-};"C";"F" **\dir{-};"E";"F"
**\dir{-};"F";"D" **\dir{-}; "F";"G" **\dir{-};"G";"B" **\dir{-};"G";"H" **\dir{--};"H";"A" **\dir{--};"H";"E"
**\dir{--};"G";"A" **\dir{--};"F";"B" **\dir{-};"A";"D" **\dir{-};"F";"A" **\dir{--};"E";"A" **\dir{--};"F";"H"
**\dir{--};
\endxy
\]
There are six 3-simplex. And so the invariant is the following :
\begin{align*}
TV_{\Z_n,\alpha}(S^1\times S^1 \times
S^1)&=\frac{1}{n}\sum_{\begin{array}{c} \scriptstyle g,h,k \in
\Z_n, \\
\scriptstyle [g,h]=[g,k]=[h,k]=1
\end{array}}\frac{\alpha(g,h,k)\alpha(h,k,g)\alpha(k,g,h)}{\alpha(g,k,h)\alpha(h,g,k)\alpha(k,h,g)}\\
&=\frac{1}{n}\sum_{g,h,k \in
\Z_n}\frac{\alpha(g,h,k)\alpha(h,k,g)\alpha(k,g,h)}{\alpha(g,k,h)\alpha(h,g,k)\alpha(k,h,g)}
\end{align*}
For every finite group $G$ and for every 3-cocycle $\alpha \in
Z^3(G,k^{\star})$, we can define :
$$\beta(g,h,k)=\frac{\alpha(g,h,k)\alpha(h,k,g)\alpha(k,g,h)}{\alpha(g,k,h)\alpha(h,g,k)\alpha(k,h,g)},$$
it verifies some properties :
\begin{lem}\label{beta}
\begin{itemize}
\item[]
\item[(i)]For every finite group $G$ and for every $g,h,k \in G$, there is an action of $S_3$ over $\beta$ by permutating the terms and if $\sigma\in S_3$ :
$\sigma.\beta(g,h,k)=\beta(g,h,k)^{\epsilon(\sigma)}$, $\epsilon$ is
the signature.
\item[(ii)]if $G$ is abelian and $\alpha$ is \ref{alphabrac} then for all $g,h,k \in G$ : $\beta(g,h,k)=1$,
\end{itemize}
\end{lem}
{Proof (i) : } It is straightforward from the definition of $\beta$,
we give only the calculation for the permutation $(12)$ and for the
cycle $(123)$ :
\begin{align*}
\beta(h,g,k)&=\frac{\alpha(h,g,k)\alpha(g,k,h)\alpha(k,h,g)}{\alpha(h,k,g)\alpha(g,h,k)\alpha(k,g,h)} \\
&=\frac{1}{\frac{\alpha(g,h,k)\alpha(h,k,g)\alpha(k,g,h)}{\alpha(h,g,k)\alpha(g,k,h)\alpha(k,h,g)}} \\
&=\beta(g,h,k)^{-1} \\
\beta(h,k,g)&=\frac{\alpha(h,k,g)\alpha(k,g,h)\alpha(g,h,k)}{\alpha(h,g,k)\alpha(k,h,g)\alpha(g,k,h)}\\
 &=\beta(g,h,k)\qed
\end{align*}

{Proof (ii) : } \qed

Thus for $G=\Z_n$ and for $\alpha$ defined by (\ref{zn}), we have :
$$TV_{\Z_n,\alpha}(S^1\times S^1 \times S^1)=\frac{1}{n}\sum_{g,h,k \in \Z_n}1=n^2.$$
In the same way, we can compute $TV_{\Z_n,\alpha}(\Sigma_g\times
S^1)$, where $\Sigma_g$ is the closed surface of genus $g$, using
the following triangulation of $\Sigma_g$ :
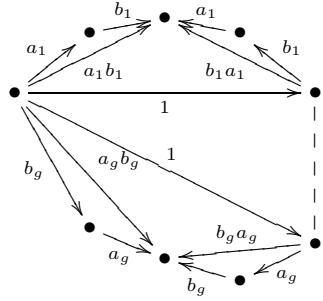
\[
\xy
(0,0)*+{\bullet}="A";(10,8)*+{\bullet}="B";(20,10)*+{\bullet}="C";(30,8)*+{\bullet}="D";(40,0)*+{\bullet}="E";(40,-20)*+{\bullet}="F";(30,-25)*+{\bullet}="G";(10,-18)*+{\bullet}="I";(20,-22)*+{\bullet}="H";
{\ar^{a_1} "A";"B"}; {\ar^{b_1} "B";"C"}; {\ar_{a_1} "D";"C"};
{\ar_{b_1} "E";"D"}; "E";"F" **\dir{--}; {\ar_{b_g}
"A";"I"};{\ar_{a_g} "I";"H"}; {\ar^{a_g} "F";"G"};{\ar^{b_g}
"G";"H"};{\ar_{a_1b_1} "A";"C"};{\ar^{b_1a_1} "E";"C"};{\ar_{1}
"A";"E"};{\ar^{a_gb_g} "A";"H"};{\ar_{b_ga_g} "F";"H"};{\ar^{1}
"A";"F"};
\endxy
\]
The edges with the same labels are identified, we have $1$ on the
edges inside the polygon by commutativity of the group $\Z_n$ and
the condition $(ii)$ of the admissible coloration. We denote $T_g$
the previous triangulation of $\Sigma_g$ : $Col(T_g)=k[\Z_n^{2g}]$.
Thus $c=(a_1,b_1,...,a_g,b_g)\in \Z_n\times ... \times \Z_n$ is a
coloring of $T_g$ We can define a triangulation of $\Sigma_g\times
S^1$ :
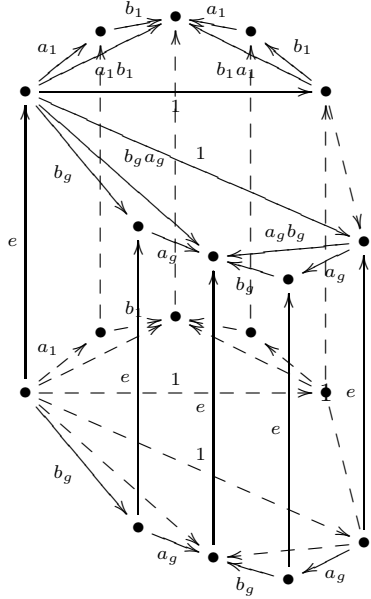
\[
\xy
(0,0)*+{\bullet}="A";(10,8)*+{\bullet}="B";(20,10)*+{\bullet}="C";(30,8)*+{\bullet}="D";(40,0)*+{\bullet}="E";(45,-20)*+{\bullet}="F";(35,-25)*+{\bullet}="G";(15,-18)*+{\bullet}="I";(25,-22)*+{\bullet}="H";(0,-40)*+{\bullet}="J";(10,-32)*+{\bullet}="K";(20,-30)*+{\bullet}="L";(30,-32)*+{\bullet}="M";(40,-40)*+{\bullet}="N";(45,-60)*+{\bullet}="O";(35,-65)*+{\bullet}="P";(15,-58)*+{\bullet}="Q";(25,-62)*+{\bullet}="R";
{\ar^{a_1} "A";"B"}; {\ar^{b_1} "B";"C"}; {\ar_{a_1} "D";"C"};
{\ar_{b_1} "E";"D"};{\ar@{-->} "E";"F"}; {\ar_{b_g}
"A";"I"};{\ar_{a_g} "I";"H"}; {\ar^{a_g} "F";"G"};{\ar^{b_g}
"G";"H"};{\ar_{a_1b_1} "A";"C"};{\ar^{b_1a_1} "E";"C"};{\ar_{1}
"A";"E"};{\ar^{b_ga_g} "A";"H"};{\ar_{a_gb_g} "F";"H"};{\ar^{1}
"A";"F"};{\ar^{e} "J";"A"};{\ar@{-->}
"K";"B"};{\ar@{-->}"L";"C"};{\ar@{-->}^{a_1} "J";"K"
};{\ar@{-->}^{b_1}
"K";"L"};{\ar@{-->}"M";"L"};{\ar@{-->}"N";"M"};{\ar@{-->}"M";"D"};{\ar@{-->}
"N";"E" };{\ar^{a_g} "O";"P"};{\ar^{b_g} "P";"R"};{\ar_{a_g}
"Q";"R"};{\ar_{b_g} "J";"Q"};"N";"O" **\dir{--};{\ar^{e}
"Q";"I"};{\ar^{e} "R";"H"};{\ar^{e} "P";"G"};{\ar^{e}
"O";"F"};{\ar@{-->}"J";"R"
};{\ar@{-->}"O";"R"};{\ar@{-->}^1"J";"O"}*{1};{\ar@{-->}"J";"L"
};{\ar@{-->}"N";"L"};{\ar@{-->}^1"J";"N"};
\endxy
\]
To avoid problems of reading this triangulation, some edges are not
coloring. We can recover the corresponding coloring thanks to the
condition $(ii)$ of coloring. So by definition of the invariant :
$TV_{\Z_n,\alpha}=\frac{1}{n}\sum_{\gamma \in
Col(T'_g)}\prod_{\Delta \in T_g^3}\alpha(\Delta,\gamma)$. The
previous triangulation can be divided into prisms :
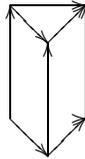
\[
\xy (0,0)*{}="A";(10,0)*{}="B";(5,-5)*{}="C";
(0,-15)*{}="D";(10,-15)*{}="E";(5,-20)*{}="F"; {\ar "D";"A"};{\ar
"A";"B"};{\ar "A";"C" };{\ar "C";"B" };{\ar "D";"F"};{\ar "F";"E"};
{\ar@{-->}|{\hole \hole} "D";"E"}; {\ar "F";"C"};{\ar "E";"B"};
\endxy
 \]
 Then we have to decompose this prism into 3-simplex. Below we give such a decomposition :
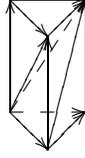
\[
\xy (0,0)*{}="A";(10,0)*{}="B";(5,-5)*{}="C";
(0,-15)*{}="D";(10,-15)*{}="E";(5,-20)*{}="F"; {\ar "D";"A"};{\ar
"A";"B"};{\ar "A";"C" };{\ar "C";"B" };{\ar "D";"F"};{\ar "F";"E"};
{\ar@{-->}"D";"E"};{\ar "F";"C"};{\ar "E";"B"};{\ar "F";"B"};{\ar
"D";"C"};{\ar@{-->} "D";"B"};
\endxy
 \]
If we assigns a coloring $\gamma$ to this prism we can compute the
term $\alpha(\Delta,\gamma)$ of the prism and by construction we
will have the invariant $TV(\Sigma_g \times I)$.
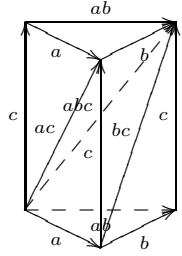
\[
\xy (0,0)*{}="A";(20,0)*{}="B";(10,-5)*{}="C";
(0,-25)*{}="D";(20,-25)*{}="E";(10,-30)*{}="F"; {\ar^{c}
"D";"A"};{\ar^{ab} "A";"B"};{\ar_{a} "A";"C" };{\ar_{b} "C";"B"
};{\ar_{a} "D";"F"};{\ar_{b} "F";"E"};
{\ar@{-->}_{ab}"D";"E"};{\ar^{c} "F";"C"};{\ar^{c}
"E";"B"};{\ar^{bc} "F";"B"};{\ar^{ac} "D";"C"};{\ar@{-->}^{abc}
"D";"B"};
\endxy
 \]
$(a,b,c)$ is a coloring of the prism, the scalar assigns to the
prism is then : $\frac{\alpha(a,b,c)\alpha(c,a,b)}{\alpha(a,c,b)}$.
In the triangulation of $\Sigma_g\times S^1$ for a given coloring,
we have prism such that the coloring is : $(a,a^{-1},c)$, where $a,c
\in \Z_n$ and so there are at least two edges whose the associated
coloring is equal to $I$. The value of the scalar associated to the
prism equip with a coloration $(a,a^{-1},c)$ is given by the
following lemma :
\begin{lem}
If $G$ is a finite group and $\gamma$ is a coloring of a 3-simplex
(oriented)$\Delta$ such that the coloring of one edge (at least) is
$1$ then for every normalized 3-cocycle in the group cohomology
$H^3(G,k^{\star})$, $\alpha(\Delta,\gamma)=1$
\end{lem}
{Proof : } If we fix $\gamma$ and $\Delta$ we have the following
figure :
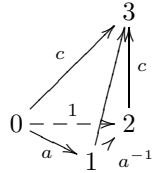
\[
\xy
(0,0)*+{0}="A";(15,0)*+{2}="B";(10,-5)*+{1}="C";(15,15)*+{3}="D";
{\ar@{-->}^{1} "A";"B"};{\ar_{a} "A";"C"};{\ar_{a^{-1}}
"C";"B"};{\ar^{c} "A";"D"};{\ar_{c} "B";"D"};{\ar "C";"D"};
\endxy
\]
The numbering of the vertices gives an orientation of $\Delta$, for
an orientation given by the numbering $0<1<2<3$ $\Delta$ will be
noted : $(0123)$. In this case we have : $\alpha(a,a^{-1},c)=1$. It
follows the definition of $H^3(G,k^{\star})$ which isomorphic to :
$H^3(BG,k^{\star})$. If we change the orientation (presevering or
not) we will have : $\alpha(a,a^{-1},c)$ (preserving the
orientation) or $\alpha(a,a^{-1},c)^{-1}$ \cite{Wakui}. \qed
 Thus the invariant is equal to :
\begin{eqnarray}\label{surfabe}
TV_{\Z_n,\alpha}(\Sigma_g\times
S^1)&=&\frac{1}{n}\sum_{(a_1,b_1,...,a_g,b_g,e)\in\Z_n}\beta(a_1,b_1,e)\times
....\times \beta(a_g,b_g,e) \nonumber\\
&=&\frac{1}{n}\sum_{(a_1,b_1,...,a_g,b_g,e)\in\Z_n}1 \nonumber\\
&=&n^{2g}
\end{eqnarray}
\subsection{$G=S_3$}
In the Annex B, we give a way of building 3-cocycle, and this
construction leads to the 3-cocycle given in \cite{Wakui} and
\cite{Moore} for $\Z_n$. And for $S_3$ we have
$\alpha(x,y,z)=exp(\frac{2i\pi}{4}tr(s(x))(tr(s(y)s(z)s(yz)^{-1})))$.
We denote $\z(x)=\{g\in S_3 \mid gxg^{-1}=x\}$ the center of $x$ we
have :
\begin{align*}
TV_{S_3,\alpha}(S^1\times S^1 \times S^1)&=\frac{1}{6}\sum_{g,h,k
\in S_3,
[g,h]=[g,k]=[h,k]=1}\frac{\alpha(g,h,k)\alpha(h,k,g)\alpha(k,g,h)}{\alpha(g,k,h)\alpha(h,g,k)\alpha(k,h,g)} \\
&=\frac{1}{6}(\sum_{h\in S_3}\sharp Z(h)+\sum_{h\not=1}(\sharp Z(h))^{2}) \\
&=\frac{1}{\sharp S_3}\sum_{g,h,k \in S_3, [g,h]=[g,k]=[h,k]=1}1 \\
&=\frac{\sharp Col(T_0)}{\sharp S_3} \\
&=\sharp\frac{Hom(\p (S^1\times S^1\times S^1),S_3)}{conj}
\end{align*}
In fact, for every 3-cocycle $\alpha \in Z^3(S_3,k^{\star})$ and
thanks to (\ref{beta}) :
\begin{align*}
\beta(g,g,k)&=\frac{\alpha(g,g,k)\alpha(g,k,g)\alpha(k,g,g)}{\alpha(g,k,g)\alpha(g,g,k)\alpha(k,g,g)}\\
&=1\\
\beta(g,h,h)&=\frac{\alpha(g,h,h)\alpha(h,h,g)\alpha(h,g,h)}{\alpha(g,h,h)\alpha(h,g,h)\alpha(h,h,g)} \\
&=1\\
\beta(g,h,g)&=\frac{\alpha(g,h,g)\alpha(h,g,g)\alpha(g,g,h)}{\alpha(g,g,h)\alpha(h,g,g)\alpha(g,h,g)} \\
&=1
\end{align*}
Moreover if $\alpha$ is normalized then $\beta$ becomes normalized.
Thus for every $\alpha \in Z^3(G,k^{\star})$ :
\begin{align*}
TV_{S_3}(S^1\times S^1\times S^1)&=\frac{1}{6}\sum_{g,h,k \in S_3, [g,h]=[g,k]=[h,k]=1}\beta(g,h,k) \\
&=\frac{1}{6}\sum_{g\in S_3}\sum_{h\in \z(g)}\sum_{k\in Z(g)\cap Z(h)}\beta(g,h,k) \\
&=\frac{1}{6}(\sum_{h\in S_3}\sum_{k\in Z(h)}\beta(1,h,k)+\sum_{g\not=1}(\sharp \z(g))^2) \\
&=\frac{1}{6}(\sum_{h\in S_3}\sharp \z(h)+\sum_{g\not=1}(\sharp \z(g))^2) \\
&=\sharp\frac{Hom(\p (S^1\times S^1\times S^1),S_3)}{conj}
\end{align*}
\subsection{Examples of TQFTs}
\subsubsection{$\alpha=1$}
Let $\Sigma$ a closed and connected surface and $T$ a triangulation
of $\Sigma$, for all $c \in Col(T)$:
\begin{align*}
V(\Sigma\times I,T)(c)&=\sum_{c' \in Col(T)}TV^n((\Sigma,T)\times I)_{c,c'}c' \\
&=\sum_{c'\cong c}\frac{\sharp stab(c)}{\sharp \lc^{ n_0}}c' \\
&=\frac{1}{\sharp \m{O}_c}\sum_{c' \cong c}c'
\end{align*}
Thus : $\m{V}(\Sigma)=\frac{k[Col(T)]}{\cong}\simeq
\frac{Hom(\p(\Sigma),\lc)}{conj}$.

\subsubsection{$G=\Z_n$, and $\alpha$ is given by (\ref{alpha})}
The vector space $\m{V}(\Sigma)$ doesn't depend on the choice of the
triangulation. Thus we can consider the following triangulation :
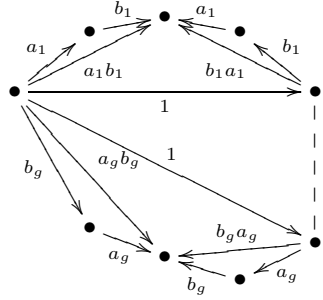
\[
\xy
(0,0)*+{\bullet}="A";(10,8)*+{\bullet}="B";(20,10)*+{\bullet}="C";(30,8)*+{\bullet}="D";(40,0)*+{\bullet}="E";(40,-20)*+{\bullet}="F";(30,-25)*+{\bullet}="G";(10,-18)*+{\bullet}="I";(20,-22)*+{\bullet}="H";
{\ar^{a_1} "A";"B"}; {\ar^{b_1} "B";"C"}; {\ar_{a_1} "D";"C"};
{\ar_{b_1} "E";"D"}; "E";"F" **\dir{--}; {\ar_{b_g}
"A";"I"};{\ar_{a_g} "I";"H"}; {\ar^{a_g} "F";"G"};{\ar^{b_g}
"G";"H"};{\ar_{a_1b_1} "A";"C"};{\ar^{b_1a_1} "E";"C"};{\ar_{1}
"A";"E"};{\ar^{a_gb_g} "A";"H"};{\ar_{b_ga_g} "F";"H"};{\ar^{1}
"A";"F"};
\endxy
\] where $(a_1,b_1,...,a_g,b_g)$ is an admissible coloring. Thus
for the cylinder $\Sigma_g\times I$ we have the following
triangulation :
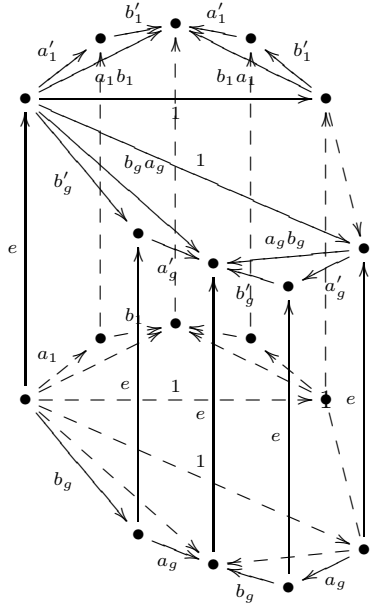
\[
\xy
(0,0)*+{\bullet}="A";(10,8)*+{\bullet}="B";(20,10)*+{\bullet}="C";(30,8)*+{\bullet}="D";(40,0)*+{\bullet}="E";(45,-20)*+{\bullet}="F";(35,-25)*+{\bullet}="G";(15,-18)*+{\bullet}="I";(25,-22)*+{\bullet}="H";(0,-40)*+{\bullet}="J";(10,-32)*+{\bullet}="K";(20,-30)*+{\bullet}="L";(30,-32)*+{\bullet}="M";(40,-40)*+{\bullet}="N";(45,-60)*+{\bullet}="O";(35,-65)*+{\bullet}="P";(15,-58)*+{\bullet}="Q";(25,-62)*+{\bullet}="R";
{\ar^{a'_1} "A";"B"}; {\ar^{b'_1} "B";"C"}; {\ar_{a'_1} "D";"C"};
{\ar_{b'_1} "E";"D"};{\ar@{-->} "E";"F"}; {\ar_{b'_g}
"A";"I"};{\ar_{a'_g} "I";"H"}; {\ar^{a'_g} "F";"G"};{\ar^{b'_g}
"G";"H"};{\ar_{a_1b_1} "A";"C"};{\ar^{b_1a_1} "E";"C"};{\ar_{1}
"A";"E"};{\ar^{b_ga_g} "A";"H"};{\ar_{a_gb_g} "F";"H"};{\ar^{1}
"A";"F"};{\ar^{e} "J";"A"};{\ar@{-->}
"K";"B"};{\ar@{-->}"L";"C"};{\ar@{-->}^{a_1} "J";"K"
};{\ar@{-->}^{b_1}
"K";"L"};{\ar@{-->}"M";"L"};{\ar@{-->}"N";"M"};{\ar@{-->}"M";"D"};{\ar@{-->}
"N";"E" };{\ar^{a_g} "O";"P"};{\ar^{b_g} "P";"R"};{\ar_{a_g}
"Q";"R"};{\ar_{b_g} "J";"Q"};"N";"O" **\dir{--};{\ar^{e}
"Q";"I"};{\ar^{e} "R";"H"};{\ar^{e} "P";"G"};{\ar^{e}
"O";"F"};{\ar@{-->}"J";"R"
};{\ar@{-->}"O";"R"};{\ar@{-->}^1"J";"O"}*{1};{\ar@{-->}"J";"L"
};{\ar@{-->}"N";"L"};{\ar@{-->}^1"J";"N"};
\endxy
\]
where $(a_1,b_1,...,a_g,b_g)$ is an admissible coloring of the
inward surface, $(a'_1,b'_1,...,a'_g,b'_g)$ is an admissible
coloring of the outward surface and $e$ is a gauge. In this case we
have :
\begin{align*}
V_{\Z_n,\alpha}(\Sigma_g \times I,T)(c)&= \sum_{c'\in
Col(T)}TV_i(\Sigma_g \times I)_{c,c'}c' \\
&= \sum_{c'\cong c}n^{-1}n c' \\
&= c
\end{align*}
Thus $\m{V}(\Sigma_g)=k[Col(T)]=k[\lc^{2g}]$
\appendix
\section{Some computation of $H^3(g,k^{\star})$}
Let $G$ a finite group and $G'$ a subgroup such that :
$$
G' \stackrel{p}{\rightarrow} G \to \{1 \}
$$
and $s$ is a section of $p$. Let $A$ an abelian group and $G$ acts
trivially on $A$, then for every application : $<\mbox{ }>\mbox{ :
}G'\times G'\rightarrow A$ which verify :
\begin{itemize}
\item[(i)]<xy,z>=<x,z>+<y,z>
\item[(ii)]<x,yz>=<x,y>+<x,z>
\end{itemize}
we can define the following application :
\begin{align}\label{alphabrac}
\alpha \mbox{ : }G\times G\times G & \rightarrow A \nonumber\\
(x,y,z) & \mapsto \alpha(x,y,z)=<s(x),s(y)s(z)s(yz)^{-1}>
\end{align}
\begin{prop}
If $<ker(p),ker(p)>=0$ then $\alpha \in Z^3(G,A)$.
\end{prop}
{Proof : } For all $x,y,z,t\in G$ we have :
\begin{align*}
&\delta(\alpha)(x,y,z,t)=\alpha(y,z,t)-\alpha(xy,z,t)+\alpha(x,yz,t)-\alpha(x,y,zt)+\alpha(x,y,z) \\
&=<s(y),s(z)s(t)s(zt)^{-1}>-<s(xy),s(z)s(t)s(zt)^{-1}>+<s(x),s(yz)s(t)s(yzt)^{-1}>\\
&-<s(x),s(y)s(zt)s(yzt)^{-1}>+<s(x),s(y)s(z)s(yz)^{-1}>\\
&=<s(y)s(xy)^{-1},s(z)s(t)s(zt)^{-1}>+<s(x),s(z)>+<s(x),s(t)>-<s(x),s(zt)> \\
&=<s(y),s(z)s(t)s(zt)^{-1}>+<s(x),s(z)s(t)s(zt)^{-1}>\\
&=<s(x)s(y)s(xy)^{-1},s(z)s(t)s(zt)^{-1}> \\
&=0 \qed
\end{align*}

\subsection{Example}
\subsubsection{$G=\Z_n$}
If we consider the following application :
\begin{align*}
\Z & \rightarrow \Z_n \\
x & \mapsto \overline{x}
\end{align*}
$s \mbox{ : }\Z_n \rightarrow$ which assigns for all $\overline{x}
\in \Z_n$ its representative element in $\{0,...,n-1 \}$ is a
section of $p$. We define :
\begin{align*}
<\mbox{ }>\Z\times \Z &\rightarrow \mathbb{C} \\
(x,y) & \mapsto exp(\frac{2i\pi}{n^2}s(x)s(y)),
\end{align*}
$<\mbox{ }>$ verifies $(i)$ and $(ii)$ and if $x \in ker(p)$ then
$s(x)=0$ thus if $x$ and $y$ are elements of $ker(p)$ then
$<x,y>=0$. It defines a 3-cocycle :
\begin{equation}\label{zn}
\alpha(x,y,z)=exp(\frac{2i\pi}{n^2}s(x)(s(y)+s(z)-s(y+z))),
\end{equation}
we recover the 3-cocycle defines in  \cite{Wakui}.\\
\subsubsection{$G=S_n$}
\begin{align*}
p \mbox{ : }B_n&\rightarrow S_n \\
\tilde{\sigma} & \mapsto \sigma
\end{align*}
and the section $s$ which assigns for all permutation in $S_n$ an
elementary braid
\begin{align*}
< \mbox{ }> \mbox{ : }B_n \times B_n &\rightarrow \mathbb{C} \\
(\sigma,\tau) & \mapsto exp(\frac{2i\pi}{4}tr(x)tr(y))
\end{align*}
with $tr(x)=\sharp(\mbox{positive crossing of
x})-\sharp(\mbox{negative crossing of x})$.
\begin{lem}
\begin{enumerate}
\item If $x\in ker(p)$ then $tr(x)$ is an even number. \\
\item $tr(xy)=tr(x)+tr(y)$.
\end{enumerate}
\end{lem}
{Proof : } We recall that $ker(p)=P_n$ the pure braids group. The
first assertion is then a consequence of the presentation of $P_n$
due to
Markov (see \cite{Markov}). \\
The last assertion is a consequence of the braid relations which preserve the number of signed crossings.\\
Thus $<\mbox{ }>$ defines a 3-cocycle $\alpha$ on $S_n$.\qed

\begin{ack}
The author would like to thank his advisor Alain Brugui\`eres for
his useful comments and constructive remarks. \\
\end{ack}
\bibliographystyle{plain}
\bibliography{biblio}
\end{document}